\documentclass{article}

           \usepackage[latin1]{inputenc}
           \usepackage[T1]{fontenc}
           \usepackage{amsfonts}
           \usepackage{amssymb}
\usepackage{amsmath}

\usepackage{graphicx}
\usepackage{psfrag}
\usepackage[dvips]{epsfig}

\newcommand{\N}{\mathbb{N}}
\newcommand{\R}{\mathbb{R}}

\newcommand{\ds}{\displaystyle}
\newtheorem{theorem}{Th\'eor\`eme}[section]

\newtheorem{proposition}[theorem]{Proposition}
\newtheorem{definition}[theorem]{D\'efinition}

\usepackage{fancyhdr}

\title{A proof of the invariance of the contact angle in Electrowetting}
\author{Claire Scheid \and Patrick Witomski}

\date{}

\begin{document}
\renewcommand{\labelitemi}{\textbullet}
\maketitle
\begin{center}
{\it LJK, Grenoble, France}
\end{center}
\begin{abstract}
We prove the invariance of the contact angle in liquid-solid wetting phenomena : an electrified droplet is spreading on a solid surface. The drop is minimizing its energy. We express the differential of this energy with respect to the shape of the drop and deduce necessary conditions for optimality . By variational method, using well-chosen test functions, we obtain the main result about the contact angle between the drop and the solid.

MSC-class : 49K10,49K20,49K30.

\end{abstract}
\maketitle
\section*{Introduction}\label{intro}
 In this article we give a proof of the invariance of the contact angle in a solid-liquid wetting phenomenon. To model this experience, described below, we use an energy minimization with respect to the shape and we obtain the main result by variational techniques.
The phenomenon we are dealing with is named electrowetting.

Electrowetting is a technique allowing to modify the affinity
between a solid and a liquid, by the introduction of an electric
field. This phenomenon has been discovered in the $XIX ^{th}$
century by Gabriel Lippmann on a mercury/electrolyte solution system (\cite{L}).
It has been then studied on more general systems (\cite{B,QB}). We can find today
many industrial applications. Particulary in optic with variable focal
lenses (see the web site of the society of B.Berge : http://www.varioptic.com/en/
and \cite{BP}), and pixels for electronic paper (\cite{pixel,pixel2}).
It is also use in microfluidics, for example in chip design (\cite{PP}) and has biomedical applications (\cite{AHT})

Consider a liquid drop on a polymer film. The equilibrium of the drop
results on superficial tension forces and gravitational force. The
drop shape is a quasi-spherical calotte with a contact angle between
water and solid which is given by Young's angle (1805)\cite{F}. We applied then a
constant voltage $\phi_0$ between this drop and an electrode placed
under the insulating polymer. The charged drop and the electrode create a capacitor.
A simple model considering the system as a plane capacitor predicts
the total spreading of the drop when the voltage becomes higher. However
physical experiences show a locking of the contact angle when the drop
is bound to a voltage higher than a critical value.

Many explanations have been proposed for this saturation of the angle.
Most of them give a crucial role to the divergence of the electric
field in the vicinity of the triple line (solid-liquid-gas interface).The understanding of this phenomenon is slowed down by the misunderstanding
of the geometry of the drop near the triple line (or wetting line).

We prove in this paper that the contact angle is independent of the applied
potential and that the value equals Young's angle, as observed in \cite{allemands}.\\

The problem will be described in the first part and the equations of
the model established.

In the second part we will recall some shape optimization results adapted to
the model.

In the third part we will establish necessary conditions for optimality which
will be exploited in a fourth part during the calculation of the contact
angle value.
\section{Description of the problem}\label{modele}
 \indent We study the evolution of the shape of a droplet located on a
substrate and subjected to a constant voltage. A 3D model (as found in \cite{SB})
allows us to examine axisymmetric case and envisage developments to non
symmetric shape for high voltages.
\subsection{Hypotheses}

Assumptions : \\

i)The applied electrical potential $\phi_0$ is continuous\\

ii)The liquid drop is a perfect conductor.\\

iii)Electrostatics effects are negligible far away from the drop.\\
%

\subsection{Experimental device}

We consider an orthonormal basis in $\R^3$. The top side of the
polymer film, where the droplet is posed, is the plane $(Oxy)$.
We use the indices L,S and G to refer to liquid, solid and gas
domains. A couple of indices LS,LG,... relates to a liquid-solid, liquid-gas...interaction

$\Omega_e$ is the bounded domain where the experimental device takes place  and
where calculations are directed.

\begin{center}
\label{figure 2}

\scalebox{0.5}{\includegraphics{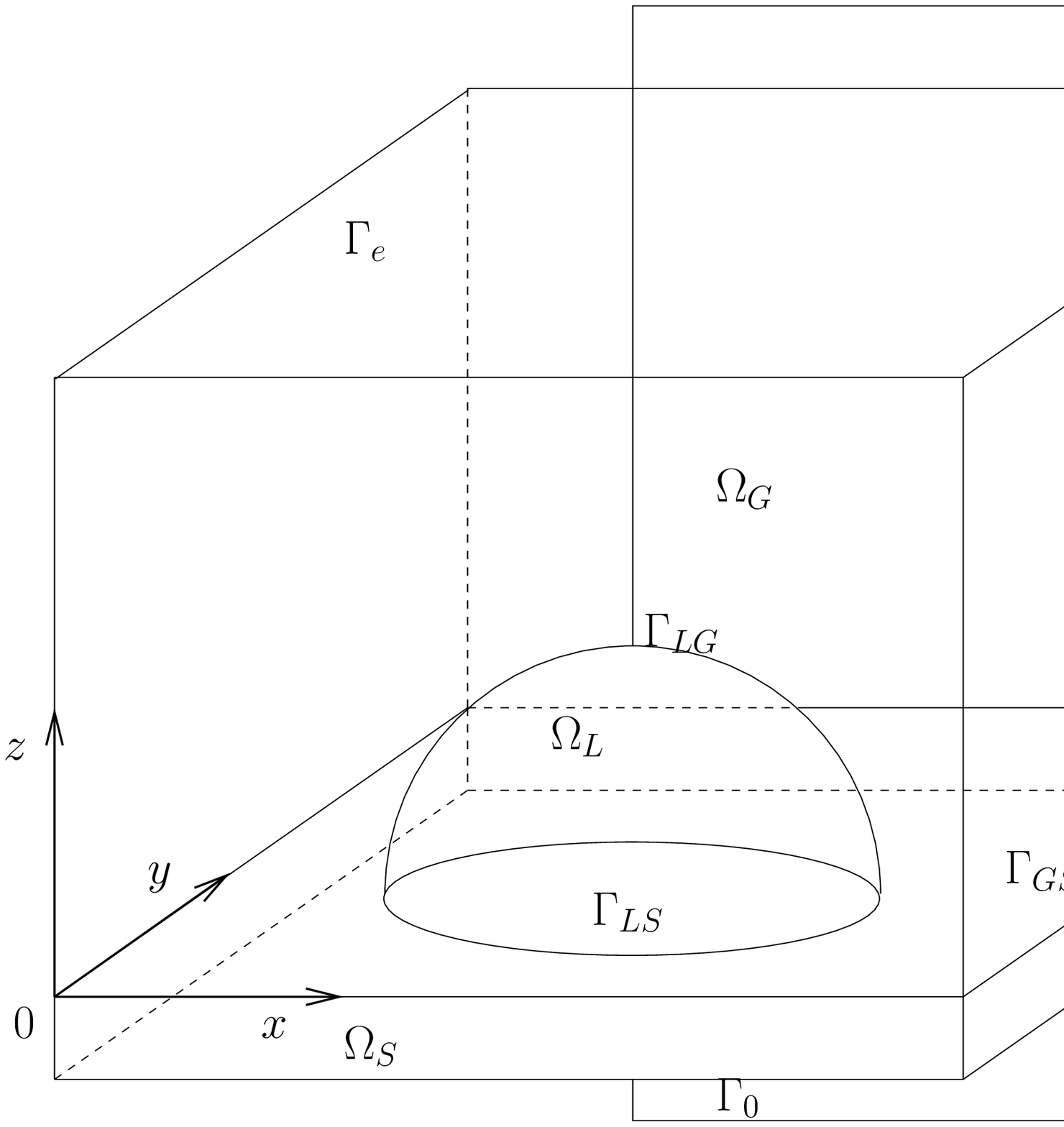}}
\end{center}

We denote by : 

$\Omega _L$ the liquid domain and $\Gamma _L$ its boundary.

$\Gamma _{LS}$ the liquid-solid interface.

$\Gamma _{LG}$ the liquid-gas interface.

$\Omega _G$ the gas domain and $\Gamma _G$ its boundary.

$\Omega _S$ the solid domain and $\Gamma _S$ its boundary.

$\Gamma _0$ the boundary of $\Omega_e$ where the counter electrode is applied.

$\Gamma _{eG}$ and $\Gamma _{eS}$ the other external boundaries of $\Omega_e$

$\Omega =\Omega _e\setminus\overline{\Omega} _L$

$\varepsilon _ G$, $\varepsilon _ S$ and $\varepsilon _ L$
permittivities of $\Omega _G$, $\Omega _S$ and $\Omega _L$, respectively.

$\overrightarrow{N_i}$, $i$ equals to $L,G,S,LS,LG,GS$...,
normals to the surfaces $\Gamma _i$, respectively.

\subsection{The electrostatic model}\label{potentiel}
 
\indent The application of an electrical potential $\phi_0$ between the
counter electrode $\Gamma_0$ and the drop $\Omega _L$ creates an
electrical potential in the entire space. The drop is supposed to be
perfectly conductive and the potential is also constant in $\Omega
_L$. But the charges distribution is not constant : it depends on the
shape of the drop i.e. on $\Omega_e\setminus\overline{\Omega}_L$.

$\phi ^{\Omega}$ is the solution of a system of partial differential
equations : 
$$\begin{array}{rl}
\textrm{div}(\varepsilon_i\nabla \phi_i^{\Omega})=0 & \mbox{ in }\Omega_i, \quad
i=G,S.\\
\phi_G^{\Omega}=\phi_0 & \mbox{ on } \Gamma_{LG}\\
\phi_S^{\Omega}=\phi_0 & \mbox{ on } \Gamma_{LS}\\
\phi_S^{\Omega}=0 & \mbox{ on } \Gamma_{0}\\
\end{array}$$

\mbox{ At the solid-gas interface, we have the following transmission
relations : } 
$$\begin{array}{cl}
\phi_G^{\Omega}=\phi_S^{\Omega} & \mbox{ on } \Gamma_{SG}\\
\varepsilon_G \nabla \phi_G^{\Omega}. \overrightarrow{N_G}=-\varepsilon_S \nabla \phi_S^{\Omega}.\overrightarrow{N_S} & \mbox{ on } \Gamma_{SG}\\
\end{array}$$

\mbox{ On the artificial boundary, we impose : }

$$\begin{array}{rl}
\varepsilon_i\nabla \phi_i^{\Omega}. \overrightarrow{N_i}=0 & \mbox{ on }\Gamma_{ei}, \quad
i=G,S.\\
\end{array}$$

This system which gives the potential can be rewritten in a weaker
form. We set :

$$H_0(\Omega)=\{\psi\in H^1(\Omega);\psi=0\mbox{ on } \Gamma _0\cup \Gamma _{LS} \cup\Gamma _{LG}\}$$

$$H(\Omega)=\{\psi\in H^1(\Omega);\psi=0\mbox{ on } \Gamma _0 \mbox{, }\psi=\phi_0\mbox{ on }\Gamma _{LS}\cup \Gamma _{LG}\}$$

Denoting $\varepsilon$ the map defined on $\Omega$ by : 

$$\varepsilon=\left\{
\begin{array}{rl}
\varepsilon _G &\mbox{ on } \Omega _G\\
\varepsilon _S &\mbox{ on } \Omega _S\\
\end{array}
\right .
$$

The weak formulation is :  

$$(FV)_{\Omega}\left\{\begin{array}{ll}
\mbox{find } \phi^{\Omega}\in H(\Omega) \mbox{ such that }\\
\\ 
\forall \psi \in H_0(\Omega), \quad \ds\int _{\Omega}\varepsilon <\nabla \phi ^{\Omega}, \nabla \psi>d\Omega=0\\
\end{array}\right .$$
where $<.,.>$ is the scalar product of $\R^3$.

In the following we denote : 
\begin{equation}\label{ff}
a_{\Omega}(\phi,\psi)=\int _{\Omega}\varepsilon <\nabla \phi, \nabla
\psi>d\Omega
\end{equation}

Thanks to Lax-Milgram's theorem we are able to prove that this problem
admits a unique solution $\phi^{\Omega}\in H(\Omega)$.

\subsection{The energy of the drop}\label{energie}
 
\indent The potential $\phi _0\geq 0$ being given, the equilibrium shape of
the drop $\Omega _L$ correspond to a minimum of the energy of the
system. We take into account the superficial tension force,
gravitational force and electrostatic force. We denote $\sigma _{LS}$,
$\sigma _{GS}$ and $\sigma _{LG}$ surface tensions relative to the
different interfaces.

For a drop $\Omega _L$ and a potential $\phi
_0$, the energy of the system is :

$$\mathcal{E}(\Omega _L, \phi _0)=
\underbrace{\rho g \int_{\Omega
_L}zd\Omega}_{\mbox{Potential energy}}+\underbrace{\int_{\Gamma_{LG}}\sigma_{LG} d\sigma+
\int_{\Gamma_{LS}}\sigma_{LS}d\sigma+\int_{\Gamma_{GS}}\sigma_{GS}
d\sigma}_{\mbox{Capillary energy}}$$
$$-\underbrace{\frac{1}{2}\int_{\Omega}\varepsilon|\nabla
\phi^{\Omega}|^2d\Omega}_{\mbox{Electrostatic energy}}$$

Up to an additive constant, we have : 

$$\mathcal{E}(\Omega _L, \phi _0)=\rho g \int_{\Omega _L}zd\Omega+\int_{\Gamma_{LG}}\sigma_{LG} d\sigma+\int_{\Gamma_{LS}}(\sigma_{LS}-\sigma_{GS})
d\sigma-\frac{1}{2}\int_{\Omega}\varepsilon|\nabla
\phi^{\Omega}|^2d\Omega$$
where the electrostatic energy is affected by a minus sign because
this energy is imposed by an exterior generator. This integral depends
on $\Omega _L$ since $\Omega=\Omega _e\setminus\overline{\Omega _L}$ and $\phi
^{\Omega}$ is the solution of the variational formulation $(FV)_{\Omega}$.\\

 \subsection{Search of the optimal shape}\label{form_opt}
 \indent The drop $\Omega _L$ of volume $V$ is submitted to a voltage $\phi _0$.
The optimal shape $\Omega _L ^*$ of the drop correspond to a minimum
of the energy : 

$$\mathcal{E}(\Omega ^*_L,\phi
_0)=\min_{\{\Omega_L ;\textrm{Volume}(\Omega _L)=V\}}\mathcal{E}(\Omega
_L, \phi_0)$$

The evaluation of the function $\mathcal{E}(\Omega _L, \phi _0)$
requires the resolution of a partial differential equation's problem on
the domain $\Omega=\Omega_e\setminus\overline{\Omega_L}$. $\Omega _S$
is fixed but $\Omega
_G$ and also $\Omega$ depends on $\Omega _L$ : It is a major
difficulty of the problem.

Obviously it is equivalent to give $\Omega _L$ or to give $\Omega$,
and in the following we take $\Omega$ as a variable.

By introducing the parameters : \\
$$\alpha=\dfrac{\rho g}{\sigma_{LG}},\quad
\mu=\dfrac{\sigma_{LS}-\sigma_{GS}}{\sigma_{LG}},\quad
\delta=\dfrac{1}{\sigma_{LG}}$$ our problem is equivalent to minimizing

$$J(\Omega)=-\alpha\int_{\Omega}zd\Omega+\mu\int_{\Gamma_{LS}}d\sigma+\int_{\Gamma_{LG}}d\sigma
-\frac{\delta}{2}\int_{\Omega}\varepsilon|\nabla \phi^{\Omega}|^2d\Omega$$

We want to minimize $J$ on a set of admissible domains obtained by
small smooth deformations of a reference domain. This point will be
clarified on paragraph 2.

The goal is then to determine a necessary condition for optimality for a
domain $\Omega$. That's why we are going to give a sense to the
derivative of the functional $J$. The main difficulty is that we want
to derive on a set of domains which has not the usual stucture
required to define a derivative in the ordinary sense.
 
 \section{Necessary conditions for optimality }\label{cn_opt}
 \indent To obtain a necessary condition for optimality, let us make precise the class of
domain on which we are minimizing.

The theory we use is detailed in \cite{MS}, \cite{HP} and \cite{S}.

\subsection{Admissible domains}

\indent We need to give a sense to integral formulation on domain or on
boundary of domains. We will work with open sets $\Omega$ with
Lipschitz boundary.
For this reason, we choose to take deformations of open sets with
Lipschitz boundary, obtained by sufficient smooth maps, in
order to keep the lipschitzian feature of domains (\cite{S}).
\begin{itemize}
	
\item We denote $$\mathcal{C}^1(\overline{\Omega},\R^3):=\left\{U_{/\overline{\Omega}} ; U\in\mathcal{C}^1(\R^3,\R^3)\right \}$$
	with the infinity or sup norm $$||U||_{\infty}=\sup
	_{x\in\overline{\Omega}} |U(x)|+\sup
	_{x\in\overline{\Omega}} |DU(x)|$$
where $DU$ is the differential of $U$\\
	
As we only consider bounded domains, for all $U$ in $\mathcal{C}^1(\overline{\Omega},\R^3)$, we have $U$ and $DU$ uniformly bounded on $\bar\Omega$.

\item
We introduce the set of admissible displacements : 
$$\mathcal{U}(\overline{\Omega},\R^3)=\left\{U\in\mathcal{C}^1(\overline{\Omega},\R^3)
; ||U||_{\infty}<1 , {U_z}_{/\Omega _S}\equiv 0,\mbox{ and }U_{/\Gamma
_e}\equiv 0\right\}$$

Finally for $U\in\mathcal{C}^1(\overline{\Omega},\R^3)$, 
$\Omega+U$ denote the set $(Id+U)(\Omega)$	
\\ \item
Let $\Omega ^0$ be a fixed reference domain of the type described in
paragraph 1.
Define $$\mathcal{D}_{ad}:=\left\{\Omega^0+U, U\in
\mathcal{U}(\overline{\Omega ^0},\R^3)\right\}$$

We are searching a necessary condition for optimality for the solution of
the following problem : \\
$$(P)\left\{\begin{array}{l}
\mbox{Find }
\Omega^*\in\mathcal{D}_{ad}\mbox{ such that }\\
\\
\ds J(\Omega
^*)=\min_{\Omega\in\mathcal{D}_{ad};C(\Omega)=0}J(\Omega)\\
\end{array}\right .$$\\
where $C(\Omega)$ is the volume constraint, more precisely : \\
$C(\Omega)=\mbox{Vol}(\Omega_e\setminus\overline{\Omega})-V$, recalling that $\Omega _L=\Omega
_e\setminus\overline{\Omega}$.\\

We are now going to give a sense to the differentiation of the
function $J$ comparatively to a domain $\Omega$ with the concept of
shape derivative. It is a classical notion which can be find in detail
in \cite{MS}. Here a weaker version is given, which is still sufficient
for our problem.

   \item Directional derivative.
\begin{definition}\label{derivee_forme}
A function $J$ defined on $\mathcal{D}_{ad}$ and with values in $\R$
has a directional derivative at a point $\Omega$ of
$\mathcal{D}_{ad}$ in a direction $U\in\mathcal{U}(\bar\Omega,\R^3)$ if the function 
$\begin{array}{cccc}
J^* : &V&\mapsto J((I+V)(\Omega))
\end{array}$
 defined on $\mathcal{U}(\bar\Omega,\R^3)$ with values in $\R$ has a
 directional derivative at the point 0 in the direction $U$ (in the
 usual sense in 
$\mathcal{C}^1(\R^3,\R^3)$). The directional derivative of $J$ at
 $\Omega$ in the direction $U$ is denoted : 
$DJ(\Omega).U:=DJ^*(0).U$.
  
\end{definition}
\end{itemize}

We are now able to write a necessary condition for optimality for $J$.

For $\Omega\in\mathcal{D}_{ad}$, and $\lambda\in\R$, we denote
$\mathcal{L}(\Omega,\lambda)=J(\Omega)-\lambda C(\Omega)$ the lagrangian of our problem of optimisation under constraint.
We are going to find a saddle point of $\mathcal{L}$ (\cite{FG}). We are searching a necessary condition for optimality  for a couple $(\Omega ^*,\lambda
^*)$ to be a saddle point of $\mathcal{L}$ with $\Omega ^*\in\mathcal{D}_{ad}$ and $\lambda \in\R$.

\subsection{Necessary condition for optimality }

\begin{proposition}\label{condition_necessaire}
If $(\Omega^*,\lambda^*)$ is a saddle point of $\mathcal{L}$, then if $J$ and $C$ admit a directional derivative at $\Omega^*$ in the direction $U\in\mathcal{U}(\overline{\Omega ^*},\R^3)$, $DJ(\Omega^*).U=\lambda^*DC(\Omega^*).U$.
\end{proposition}
Proof : 
$(\Omega^*,\lambda^*)$ saddle point of  $\mathcal{L}\Leftrightarrow\forall\Omega\in\mathcal{D}_{ad},\quad\forall\lambda\in\R$,
$$\mathcal{L}(\Omega^*,\lambda)\leq\mathcal{L}(\Omega^*,\lambda^*)\leq\mathcal{L}(\Omega,\lambda^*)$$
Let $U\in\mathcal{U}(\overline{\Omega^*},\R^3)$. 
For $t\in [-1,1]$, $tU$ is an element of  $\mathcal{U}(\overline{\Omega^*},\R^3)$. And so the set $\Omega^*+tU$ is an element of $\mathcal{D}_{ad}$ too.\\
Then we have : 
$$\mathcal{L}(\Omega^*,\lambda^*)\leq\mathcal{L}(\Omega^*+tU,\lambda^*)
,\quad\forall t\in [-1,1]$$
By making explicit the values of $\mathcal{L}$, we get : 
$$J(\Omega^*)+\lambda^*C(\Omega^*)-J(\Omega^*+tU)-\lambda^*C(\Omega^*+tU)\leq 0
,\quad\forall t\in [-1,1]$$
By taking definition's notation, we deduce that : 
$$J^*(tU)-J^*(0)+\lambda^*\left [C^*(tU)-C^*(0)\right ]\geq 0
,\quad\forall t\in [-1,1]$$
Taking $t>0$, we obtain : 
$$\frac{J^*(tU)-J^*(0)}{t}+\lambda^*\frac{C^*(tU)-C^*(0)}{t}\geq 0
,\quad\forall t\in [-1,1]$$
Let $t$ tend to 0, we obtain : 
$$DJ^*(0).U+\lambda^*DC^*(0).U\geq 0$$\\
With the same argument but with $t<0$, we finally obtain : 
$$\forall U\in\mathcal{U}(\overline{\Omega^*},\R^3),\quad
DJ^*(0).U+\lambda^*DC^*(0).U=0$$
That is to say, by definition of shape derivative :
$$\forall U\in\mathcal{U}(\overline{\Omega^*},\R^3),\quad DJ(\Omega^*).U+\lambda^*DC(\Omega^*).U=0$$

\begin{flushright}
$\blacksquare$
\end{flushright}

We are now searching the expressions of the derivative of $J$ and $C$ with respect to a domain.

 \section{Derivative of the drop's energy with respect to its shape}\label{derive_forme}

\indent Let $\Omega^*\subset\mathcal{D}_{ad}$, such that 
$$J(\Omega^*)=\min_{\Omega\in\mathcal{D}_{ad} \atop C(\Omega)=0}J(\Omega)$$
where
$$J(\Omega)= \underbrace{-\alpha\int_{\Omega}zd\Omega}_{J_{grav}}+\underbrace{\int_{\Gamma_{LG}}d\sigma}_{J_{LG}}+\underbrace{\mu\int_{\Gamma_{LS}}
d\sigma}_{J_{LS}}-\underbrace{\frac{\delta}{2}\int_{\Omega}\varepsilon|\nabla
\phi^{\Omega}|^2d\Omega}_{J_{el}}$$
and
$$C(\Omega)=\mbox{Vol}(\Omega_e\setminus\overline{\Omega})-V$$

$J$ is the sum of four terms $J_{grav}, J_{LG}, J_{LS}, J_{el}$. The first three terms of $J$ and $C$ are integrals on a surface or a domain of a fonction independant of this domain. We have classical results on the shape derivation of such terms. (see for exemple \cite{MS}). 
So we won't give more details for the derivation of this terms.

We are going to spend more time on the singular term of our problem : the electrostatic contribution $J_{el}$.

\subsection{Derivative of gravitational, capillary and volume constraint term}

\indent Let $\Omega\in\mathcal{D}_{ad}$ and $U\in\mathcal{U}(\overline{\Omega},\R^3)$.

The derivatives at the point $\Omega$
and in the direction $U$ is obtained by following Definition  \ref{derivee_forme} and the results of \cite{MS} : 
\begin{itemize}
\item
\begin{equation}\label{derivee_gravite}
DJ_{grav}(\Omega).U=DJ_{grav}^*(0).U=\alpha\int _{\Omega}U_{z} d \Omega+ \alpha\int
_{\Omega} z\textrm{div}(U)d\Omega
\end{equation}

\item
In the same way, we express the differential of $C$ which looks like the term $J_{grav}$.
\begin{equation}\label{derivee_contrainte}
DC(\Omega).U=\int_{\Omega}\textrm{div}(U)d\Omega
\end{equation}

\item
The derivative of the terms $J_{LG}$ and $J_{LS}$ are obtained likewise using a result about derivative with respect to a surface : \begin{equation}\label{derivee_bord_lg}
DJ_{LG}(\Omega).U=\int _{\Gamma
  _{LG}}\textrm{div}(U)d\sigma-\int _{\Gamma
  _{LG}}<\overrightarrow {N_{LG}},^tDU\overrightarrow {N_{LG}}>d
  \sigma
\end{equation}

\begin{equation}\label{derivee_bord_ls}
DJ_{LS}(\Omega).U=\mu\int _{\Gamma
  _{LS}}\textrm{div}(U)d\sigma-\mu\int _{\Gamma
  _{LS}}<\overrightarrow {N_{LS}},^tDU\overrightarrow {N_{LS}}>d\sigma
\end{equation}

where $^tDU$ is the transposition of the Jacobian of $U$

\end{itemize}

\subsection{Derivative of the electrostatic energy}

\indent We are studying much more in detail the electrostatic term $J_{el}$ (\cite{MW}).

Let $V\in\mathcal{U}(\overline{\Omega},\R^3)$, by Definition (\ref{derivee_forme})
$$J_{el}^*(V)=\frac{\delta}{2}\int _{\Omega+V}\varepsilon \left
|\nabla \phi ^{\Omega+V} \right |^2 d\Omega$$

To lighten the notation, let us pose $F:=Id+V$.

As $V\in\mathcal{U}(\overline{\Omega},\R^3)$, $F$ is invertible.

\begin{itemize}
\item $\phi ^{F(\Omega)}\in H(F(\Omega))$ is the solution of the variational problem $(FV)_{F(\Omega)}$.

\item The map $T_0 : \psi\in H_0(F(\Omega))\mapsto \psi\circ F\in
H_0(\Omega)$ is an isomorphism from $H_0(F(\Omega))$ to $H_0(\Omega)$. Likewise we define $T$ from $H(F(\Omega))$ to
$H(\Omega)$, because $\phi_0$ is a constant ; $T$ is an isomorphism too.

We consider the two transported variational problems : 

$$(FV)_{F(\Omega)}\left\{\begin{array}{ll}
\mbox{ Find } \phi^{F(\Omega)}\in H(F(\Omega)) \mbox{ such that }\\
\\ 
\forall \psi \in H_0(F(\Omega)), \quad
a_{F(\Omega)}(\phi^{F(\Omega)},\psi)=\int _{F(\Omega)}\varepsilon<\nabla \phi^{F(\Omega)},
\nabla \psi>d\Omega=0\\
\end{array}\right .
$$

$$(FV)^T_{\Omega}\left\{\begin{array}{ll}
\mbox{ Find } {\phi}^{F}\in H(\Omega) \mbox{ such that }\\
\\ 
\forall \psi \in H_0(\Omega), \quad \int _{\Omega}\varepsilon  <^t(DF^{-1}\circ
F)\nabla { \phi}^{F},^t(DF^{-1}\circ
F)\nabla  \psi>|\det(DF)|d\Omega=0\\
\end{array}\right .
$$

Problems $(FV)_{F(\Omega)}$ and $(FV)^T_{\Omega}$ are equivalent.

In the following we note, for $ \phi\in H(\Omega)$ and $
\psi\in H_0(\Omega)$
\begin{equation}\label{a_tilde}
\tilde a (F,{ \phi}, \psi):=\int _{\Omega}\varepsilon  <^t(DF^{-1}\circ
F)\nabla { \phi},^t(DF^{-1}\circ
F)\nabla  \psi>|\det(DF)|d\Omega
\end{equation}

By unicity of solutions of the two variational problems
$${ \phi}^{F}=\phi^{F(\Omega)}\circ F$$

and so 
$$\int _{F(\Omega)}\varepsilon \left
|\nabla \phi ^{F(\Omega)} \right |^2 d\Omega=\int _{\Omega
}\varepsilon \left | ^t(DF^{-1}\circ
F)\nabla { \phi}^{F}\right | ^2  |\det(DF)|d\Omega$$

Thus, let us pose for $\phi\in H(\Omega)$, 
\begin{equation}\label{J_tilde}
\tilde{J_{el}}(F,\phi)=\frac{\delta}{2}\int _{\Omega
}\varepsilon \left | ^t(DF^{-1}\circ
F)\nabla \phi\right | ^2  |\det(DF)|d\Omega
\end{equation}
In the expression of $J_{el}^*$, variables $F$ and $\phi^F$ have been decoupled.

By the implicit function theorem we can show that $\tilde{J_{el}}$ is  $\mathcal{C}^1$ and that $V\mapsto\phi^{Id+V}$ is differentiable (\cite{S}).\\
A classical optimal control result (\cite{MS} V-10) allows us to write 
$$\forall U\in\mathcal{U}(\overline{\Omega},\R^3),\quad DJ_{el}^*(0).U=D_F\tilde{J_{el}}(Id,\phi^{\Omega}).U-D_F\tilde{a}(Id,\phi^{\Omega},p^{\Omega}).U$$
where $D_F$ denotes the partial differential with respect to the variable $F$.

$p^{\Omega}\in H_0(\Omega)$ is the adjoint state solution of the equation
\begin{equation}\label{etatadjoint}
\forall\phi\in H_0(\Omega),\quad D_{\phi}\tilde{a}(Id,\phi^{\Omega},p^{\Omega}).\phi=D_{\phi}\tilde{J_{el}}(Id,\phi^{\Omega}).\phi
\end{equation}

The differential of (\ref{a_tilde}) and (\ref{J_tilde}) with respect to 
$\phi$ when $F=Id$ gives : 
\begin{equation}\label{derivee_a_tilde}
\forall\phi\in H_0(\Omega),\quad
D_{\phi}\tilde{a}(0,\phi^{\Omega},p^{\Omega}).\phi=a_{\Omega}(p^{\Omega},\phi)
\end{equation}
and $\phi^{\Omega}$ being the solution of the variational problem $(FV)_{\Omega}$,
\begin{equation}\label{derivee_J_tilde}
D_{\phi}J_{el}(Id,\phi^{\Omega}).\phi=-\delta
a_{\Omega}(\phi^{\Omega},\phi)=0
\end{equation}
We deduce of (\ref{etatadjoint}), (\ref{derivee_a_tilde}) and (\ref{derivee_J_tilde}) that $\forall\phi\in H_0(\Omega)$,
$a_{\Omega}(p^{\Omega},\phi)=0$.

And  so as $p^{\Omega}\in H_0(\Omega)$, $$p^{\Omega}=0$$

From this we deduce,

$$\forall U\in\mathcal{U}(\overline{\Omega},\R^3),\quad DJ_{el}^*(0).U=D_F\tilde{J_{el}}(Id,\phi^{\Omega}).U$$

Now we look at an integral on a fixed domain and using the differentiation under the integral sign formulas (\cite{S}), we obtain finally
\begin{multline}\label{derivee_el}
  DJ_{el}(\Omega).U=DJ_{el}^*(0).U=-\frac{\delta}{2}\int_{\Omega}\varepsilon |\nabla {\phi}
  ^{\Omega}|^2\textrm{div}(U)d\Omega\\
+\frac{\delta}{2}\int_{\Omega}\varepsilon <(^tDU+DU)\nabla
  {\phi}^{\Omega},\nabla {\phi}^{\Omega}>d\Omega
\end{multline}
In the following, we simplify the notation by denoting $L_{grav}$,
  $L_{cont}$, $L_{LS}$, $L_{LG}$ and $L_{el}$ respectively for $DJ_{grav}$, $DC$, $DJ_{LS}$, $DJ_{LG}$ et $DJ_{el}$.

\subsection{Formulation of the necessary condition for optimality}

\noindent By proposition (\ref{condition_necessaire}), if a pair $(\Omega^*,\lambda^*)$ of $\mathcal{D}_{ad}\times\R$ is a saddle point of  $\mathcal{L}$ then 
\begin{multline}\label{formule_3d_cn}
\forall U\in\mathcal{U}(\overline{\Omega ^*},\R^3),\quad\\
L_{grav}(\Omega ^*).U
+L_{LG}(\Omega ^*).U+L_{LS}(\Omega ^*).U+L_{el}(\Omega^*).U
=\lambda^*L_{cont}(\Omega^*).U
\end{multline}

Terms used below are defined by (\ref{derivee_gravite}),
(\ref{derivee_contrainte}), (\ref{derivee_bord_lg}),
(\ref{derivee_bord_ls}) and (\ref{derivee_el}).
\end{itemize}

The formulation obtained there is verified for all 3D optimal domain of  $\mathcal{D}_{ad}$. It can be used for numerical simulations particulary by high potential $\phi_0$.

We are giving a first application to calculate the contact angle for axisymmetric shape.

 \section{Calculation of the contact angle for an axisymmetric shape}\label{calcul_angle}
 \indent This paragraph contains the proof of the main result : in the case of an axisymmetric optimal shape, the contact angle is the static Young's angle.

First, we will write necessary conditions for optimality introducing the axisymmetry in the model.

Then, we will by a judicious choice of direction of deformation calculate the value of the contact angle.

\subsection{The axisymmetric problem}

\indent Let's suppose that the domain $\Omega$ of $\R^3$ is axisymmetric. We choose to express a point of the space in cylindric coordinates.

In an orthonormal basis of $\R^2$ we denote $\omega$ the $2D$ domain associated to $\Omega$ in $3D$. We define in a similar way $\omega _L$, $\omega
_S$, $\omega _G$, $\omega_e$, $\gamma _{LS}$, $\gamma _{LG}$, $\gamma
_{GS}$, $\gamma_e$.\\
All are defined as in the figure : \\
\begin{center}
\scalebox{0.5}{\includegraphics{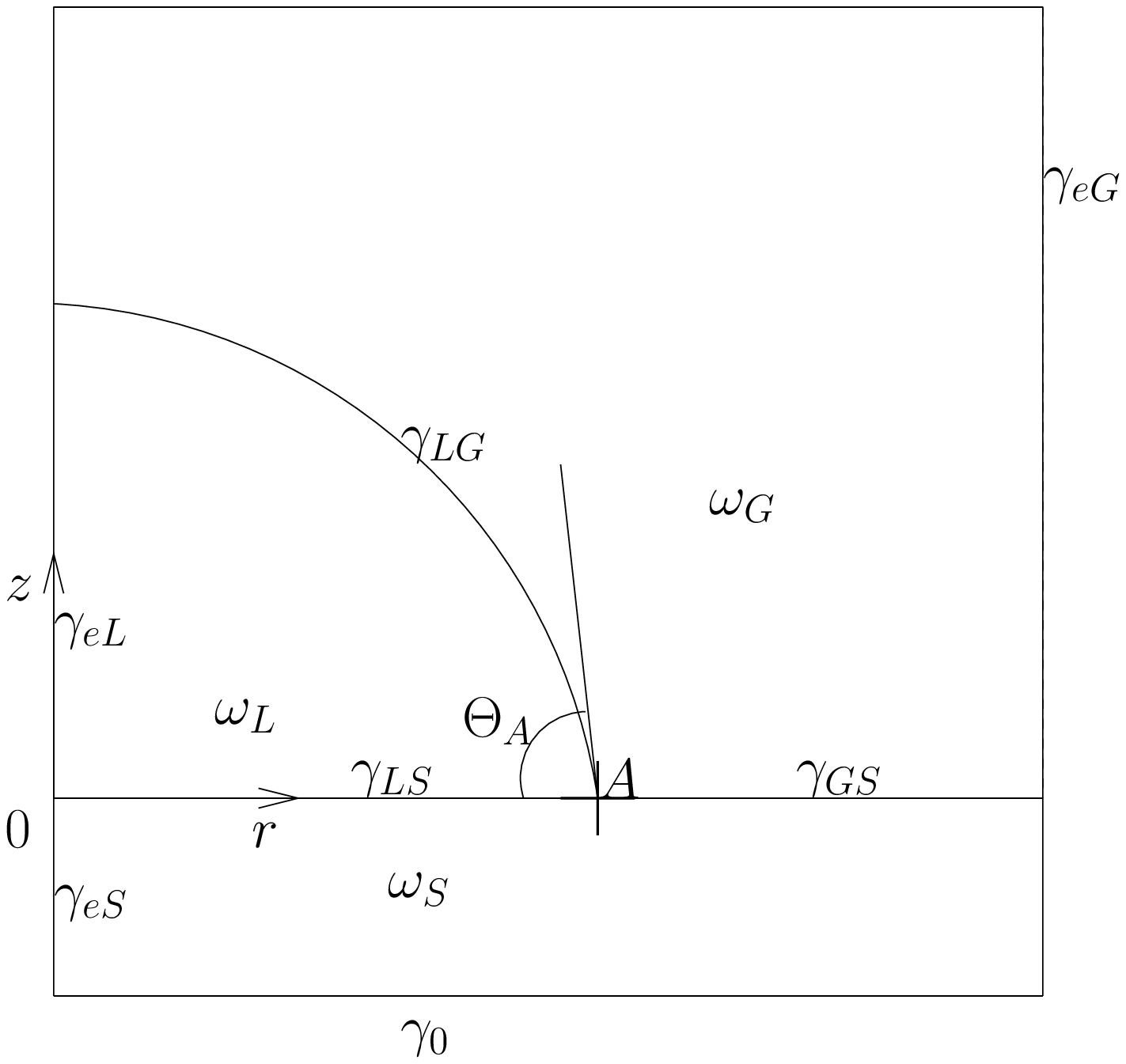}}
\end{center}
$\gamma_e$ is constituted of external liquid 
($\gamma_{eL}$), solid ($\gamma_{eS}$) and gas ($\gamma_{eG}$) boundaries.

$A$ is the contact point, and $\Theta_A$ the contact angle.
\subsubsection{Notations}

\indent For $\omega$
, we define the analogous spaces to those defined for the domain $\Omega$.

\begin{itemize}

\item

$\mathcal{C}^1(\overline{\omega},\R^2)$
\item

$\mathcal{U}(\bar\omega,\R^2)=\left\{u\in\mathcal{C}^1(\overline{\omega},\R^2)
; ||u||_{\infty}<1 , {u_z}_{/\omega _S}\equiv 0,\mbox{ et }u _{/\gamma
_e}\equiv 0\right\}$\\

Let us suppose that we are at the minimum of the energy $\Omega^*$ and consider $\omega^*$ the associated $2D$ domain.

\item
We choose a direction of deformation $U\in\mathcal{U}(\overline{\Omega ^*},\R^3)$ which is invariant by a rotation around the $z$ axis.

We are able to find
\begin{eqnarray*}
u: & \overline{\omega^*}\rightarrow &\R^2\\
   &(r,z)\mapsto &(u_r(r,z),u_z(r,z))\\
\mbox{such that} & u\in\mathcal{U}(\overline{\omega^*},\R^2) & \mbox{and if}\\ 
\tilde U: & \overline{\Omega^*}\rightarrow &\R^3\\
   &(r,\theta,z)\mapsto
   &(u_r(r,z)\cos(\theta),u_r(r,z)\sin(\theta),u_z(r,z))\\
\end{eqnarray*}
then
$$U(r\cos(\theta), r\sin(\theta),z)=\tilde U(r,\theta,z)$$

\item

In the same way, we have analogous relations and notation for the different normals to the boundaries. When there is no ambiguity on the considered surface  $\gamma$ (resp. $\Gamma$), we will denote $\overrightarrow{n}$ (resp
$\overrightarrow{N}$) for $\overrightarrow{n_{\gamma}}$ (resp
$\overrightarrow{N_{\Gamma}}$).

It exists

\begin{eqnarray*}
\overrightarrow n: & \R^2\rightarrow &\R^2\\
   &(r,z)\mapsto &(n_r(r,z),n_z(r,z))\\
\mbox{such that if} & & \\ 
\tilde N: & \R^3\rightarrow &\R^3\\
   &(r,\theta,z)\mapsto
   &(n_r(r,z)\cos(\theta),n_r(r,z)\sin(\theta),n_z(r,z))\\
\end{eqnarray*}

then 

$$\overrightarrow{N}(r\cos(\theta), r\sin(\theta),z)=\tilde N(r,\theta,z)$$

The deformation direction and the normal are also entirely determined by the maps of the plane $u$ and 
$\overrightarrow n$.

\item

We have
\begin{equation}\label{div_pol}
\textrm{div}(U)=\frac{u_r}{r}+\textrm{div}(u)
\end{equation}
 and
\begin{equation}\label{scal_pol}
<\overrightarrow N , ^tDU \overrightarrow N
>=<\overrightarrow n , ^tDu \overrightarrow n >
\end{equation}
for the normal to the different considered surfaces.

We will note $$\textrm{div} _{\gamma}(u)=\textrm{div}(u)-<\overrightarrow
n , ^tDu \overrightarrow n >$$ the surface divergence of the surface  $\gamma$.

\item

We define a $2D$ potential associated to the $3D$ potential too.

We will note
$$\begin{array}{l}
\tilde{H_0}(\omega^*)=\{ u\in H^1_{\nu}(\omega^*);u=0\mbox{ sur } \gamma^*
_0\cup \gamma^* _{LS} \cup\gamma^* _{LG}\}\\
\tilde{H}(\omega^*)=\{u\in H^1_{\nu}(\omega^*);u=0\mbox{ sur } \gamma^* _0
\mbox{, }u=\phi_0\mbox{ sur }\gamma^* _{LS}\cup \gamma^* _{LG}\} \\
\end{array}$$
where $H^1_{\nu}(\omega^*)$ is the Sobolev spaces for the measure $d\nu=rdrdz$.

We will note $L^p_{\nu}$ the analogous spaces for the case $L^p$ with the measure $\nu$.

Consider the following variational problem $$(fv)_{\omega^*}\left\{\begin{array}{l}
\mbox{Find } { \varphi}^{\omega^*}\in \tilde{H}(\omega^*)
 \mbox{ such that, }\\ 
\\
\forall v \in \tilde{H_0}(\omega^*), \quad
\int _{\omega^*}\varepsilon <\nabla \varphi^{\omega^*}, \nabla
v>r drdz=0\\
\end{array}\right .$$

 where $<.,.>$ is the scalar product of $\R^2$.
 
In the following we will denote for $\varphi\in \tilde{H}(\omega^*)$ and $v\in\tilde{H_0}(\omega^*)$,

\begin{equation}\label{a_tilde_2d}
\tilde a_{\omega^*}(\varphi,v)=\int _{\omega^*}\varepsilon r <\nabla
\varphi, \nabla v>dx
\end{equation}

To ${ \phi}^{\Omega^*}$ is also associated ${ \varphi}^{\omega^*}$ the solution of the $2D$ problem given by his weak form.

We can verify that ${ \phi}^{\Omega^*}$ is a weak solution of the problem which gives the potential on $\Omega^*$ if and only if ${ \varphi}^{\omega^*}$ is a weak solution of the associated $2D$ axisymmetric problem.

\item
Let us precise the boundary $\gamma _{LG}$ in the vicinity of the triple point.
By definition of admissible domains, $\omega^*$ is a lipschitzian open set. Let $A$ be the triple point. If we suppose that the contact angle at $A$ is in the interval $]0,\pi[$, then it exists an open ball $B_A$ centered on $A$ of radius $d_0$ and a $\mathcal{C}^1$ function $f:\R\rightarrow\R$ such that, $(r,z)\in\gamma^*_{LG}\cap B_A
\Leftrightarrow r=f(z)$ (Implicit function theorem). With such a parametrization we can express the cosine of the wetting angle at $A=(f(0),0)$ :
\begin{equation}\label{angle}
\cos(\Theta _A)=-\frac{f'(0)}{\sqrt{1+f'(0)^2}}
\end{equation}

\end{itemize}

\subsubsection{Formulation of the necessary condition for optimality in the axisymmetric case}

By adopting the previous notations, and by noting $l_{grav}$,
$l_{LG}$, $l_{LS}$, $l_{cap}$ and $l_{cont}$ the associated terms in $2D$, identity (\ref{formule_3d_cn}) becomes

\begin{multline}\label{formule_2d_cn}
\exists \lambda^*\in\R\mbox{ tel que }\forall
  u\in\mathcal{U}(\overline{\omega^*},\R^2), \quad\\
  l(\omega^*).u:=l_{grav}(\omega^*).u+l_{LG}(\omega^*).u+l_{LS}(\omega^*).u+l_{el}(\omega^*).u=\lambda^*l_{cont}(\omega^*).u
\end{multline}

We express those differents quantities in detail in the following paragraph.

\subsection{A particular choice of admissible directions}

Suppose that we are at a minimum of energy $\omega ^*$

We know that at this minimum 
$$\exists \lambda^*\in\R \mbox{ such that }\forall u
\in\mathcal{U}(\overline{\omega ^*},\R^2),
l(\omega^*).u=\lambda^*l_{cont}(\omega^*).u$$

We want to use this equality in order to extract an information on the geometry of the drop in the vicinity of the triple point $A$.
The idea is to find particular deformation's directions $u^p$, $p\in\N^*$ whose support, centered at $A$ focuses on $A$ as $p$ tends to $+\infty$.

The sequence $(u^p)_{p\in\N^*}$ which is chosen is defined by :

$$\begin{array}{lccl}
u^p : & \R^2&\rightarrow&\R^2\\
& (r,z)&\mapsto&(u_r^p,0)\\
\end{array}$$
where 
$$\begin{array}{lccl}
u^p_r : &\R^2&\rightarrow&\R\\
&(r,z)&\mapsto&\left\{\begin{array}{lc}
                \ds\exp(\frac{1}{p^2((r-f(0))^2+z^2)-1})&\mbox{ if } (r-f(0))^2+z^2<\frac{1}{p^2}\\
                0&\mbox{ else.}\\
                      \end{array}\right. 
\\
\end{array}$$

The support of $u_r^p$ is the ball of center $A$ and of radius $\displaystyle\frac{1}{p}$. It is included in the neighbourhood $B_A$ where
$\gamma_{LG}^*$ is parametrated by $f$, if $p$ is higher enough.
We have $u^p\in\mathcal{U}(\overline{\omega ^*},\R^2)$.\\

In the following paragraph, we study the limit of $l(\omega^*).u^p$ and of $l_{cont}(\omega^*).u^p$ as $p$ tends to $+\infty$.

\subsection{The limit of $l(\omega ^*).u^p$}

In the following, we denote $\ds B_p$ the open ball of center $A$ and of radius $\ds\frac{1}{p}$ and $\partial B_p$ its boundary.\\
Let us study each terms which appeared in $l(\omega^*).u^p$ and
$l_{cont}(\omega^*).u^p$

\subsubsection{The term $\quad l_{grav}(\omega^*).u^p$}

With (\ref{derivee_gravite}) and (\ref{div_pol}), we have $$\quad l_{grav}(\omega^*).u^p=\ds\alpha\int
_{\omega^*} u^p_{ z}r d rd  z+  \alpha\int
_{\omega^*} z u^p _ rd rd z + \alpha\int
_{\omega^*} z \textrm{div}(u^p)rd rd z$$

It is clear that 
\begin{equation}\label{grav_1}
\ds\int_{\omega ^*}ru^p_zdx=0\mbox{, since }u^p_z=0
\end{equation}

As $|u^p_r|\leq 1$, we have 
\begin{equation}\label{maj_zaire}
\left | \int
_{\omega^*} z u _ rd rd z \right | \leq \int
_{\omega^*\cap B_p} |z| d rd z\leq \pi\frac{1}{p^3}
\end{equation}
%

As $u^p_z= 0$, we have $\ds\int_{\omega^*}z\textrm{div}(u^p)rdrdz=\int_{\omega ^*\cap
B_p}z\frac{\partial u^p_r}{\partial r}(r,z)rdrdz$\\
By posing $R:=\sqrt{(r-f(0))^2+z^2}$
$$\frac{\partial u^p_r}{\partial
  r}(r,z)=-2p^2\frac{r-f(0)}{(p^2R^2-1)^2}\exp\left
  (\frac{1}{p^2R^2-1}\right )$$
and the upper bound $\ds\left |\frac{\partial u^p_r}{\partial
  r}(r,z)\right |\leq 2p^2\frac{R}{(p^2R^2-1)^2}\exp\left
  (\frac{1}{p^2R^2-1}\right )$\\
From the study of $\ds g :
  R>0\mapsto\frac{R}{(p^2R^2-1)^2}\exp \left (\frac{1}{p^2R^2-1}\right
  )$, we deduce
  that $\ds |g(R)|\leq \frac{C}{p}$, where $C$ is a positive real,\\
and so 
\begin{equation}\label{divergence}
\ds \left | \textrm{div}(u^p)\right |\leq 2Cp
\end{equation}

We finally obtain 
\begin{equation}\label{maj_aire2}
\left  |\int_{\omega ^*}z\textrm{div}(u^p)rdrdz\right | \leq 2Cp\left |\int_{\omega ^*\cap B _ p}
zrdrdz\right |
=2C\frac{\pi}{p^2}(f(0)+\frac{1}{p})
\end{equation}
Then with (\ref{maj_zaire}) and (\ref{maj_aire2}),

\begin{equation}\label{grav_3}
\lim _{p\to +\infty}l_{grav}(\omega ^*).u^p=0
\end{equation}

\subsubsection{The term $l_{LG}(\omega ^*).u^p$}
From (\ref{derivee_bord_lg}), (\ref{div_pol}) and (\ref{scal_pol}), we deduce 
$$l_{LG}(\omega ^*).u^p=\ds\int_{\gamma
^*_{LG}}u^p_r(s)ds+\ds\int_{\gamma ^*_{LG}}\textrm{div}
_{\gamma}(u^p)rds$$

Let $\gamma(p)$ be the intersection ordinate of $C_p$ with $\gamma
^*_{LG}$.
As $p\geq p_0$, we use the parametrization of $\gamma_{LG}^*$
by $f$.
We have 
$$\int_{\gamma
^*_{LG}}u^p_r(s)ds=\int_0^{\gamma(p)}u^p_r(f(z),z)\sqrt{1+{f'}^2(z)}dz$$
so
$$\left | \int_{\gamma ^*_{LG}}u^p_r(s)ds \right | \leq\int _0^{\gamma(p)}\sqrt{1+{f'}^2(z)}dz=l(\gamma ^*_{LG}\cap B_p)$$
and then
\begin{equation}\label{lg_1}
\lim _{p\to\infty}\int_{\gamma ^*_{LG}}u^p_r(s)ds=0.
\end{equation}

The study of the term $\int_{\gamma ^*_{LG}}\textrm{div}
_{\gamma}(u^p)rds$ which contains the surface divergence is a little bit more fastidious and contributes to the final calculation of the contact angle.

We show that the limit of this term is
$-\ds\frac{f(0)f'(0)(1+f'^2(0))^{-\frac{1}{2}}}{e}$ as $p$ tends to
$+\infty$.\\

$\ds\int_{\gamma ^*_{LG}}\textrm{div}
_{\gamma}(u^p)rds=\int_{\gamma ^*_{LG}}\textrm{div}(u^p)rds-\int_{\gamma ^*_{LG}}<^tDu^p\overrightarrow
n;\overrightarrow n>rds.$\\

As 

$$\ds^tDu^p=\left ( \begin{array}{cc}
\ds\frac{\partial u^p_r}{\partial r}&0\\
\ds\frac{\partial u^p_r}{\partial z}&0\\
\end{array}\right ) 
\mbox{ and }
\ds\overrightarrow n=\frac{1}{\sqrt{1+f'(z)^2}}\left ( \begin{array}{c}
                                         -1\\
                                         f'(z)\\
                                        \end{array}\right )$$

Therefore 
\begin{multline*}
\int_{\gamma ^*_{LG}}\textrm{div}
  _{\gamma}(u^p)rds\\
=\int_0^{\gamma(p)}\frac{\partial
    u^p_r}{\partial r}f(z)\sqrt{1+f'(z)^2}dz
-\int _0^{\gamma(p)}(\frac{\partial u^p_r}{\partial r}-f'(z)\frac{\partial u^p_r}{\partial z})\frac{1}{1+f'^2(z)}f(z)\sqrt{1+f'^2(z)}dz\\
=\int _0^{\gamma(p)}f(z)f'(z)(1+f'^2(z))^{-\frac{1}{2}}[f'(z)\frac{\partial u^p_r}{\partial r}+\frac{\partial u^p_r}{\partial z}]dz\\
=\int
_0^{\gamma(p)}f(z)f'(z)(1+f'^2(z))^{-\frac{1}{2}}\frac{\partial}{\partial
  z}[u_r^p(f(z),z)]dz\\
\end{multline*}

Let us pose $h(z):=f(z)f'(z)(1+f'^2(z))^{-\frac{1}{2}}.$\\

We are now showing that 
$$\lim_{p\to +\infty}\int_{\gamma ^*_{LG}}\textrm{div}
  _{\gamma}(u^p)rds=-\frac{1}{e}h(0)$$

As $f$ belongs to $\mathcal{C}^1([0,\gamma(p)],\R)$, $h$ is a continuous function and so particulary continuous in 0 : 
$$\forall\varepsilon>0,\quad\exists\delta(\varepsilon)>0\mbox{ such that if }|z|<\delta(\varepsilon),\mbox{ then }|h(z)-h(0)|<\varepsilon$$

We have : 
$$\left |\int _0^{\gamma(p)}h(z)\frac{\partial}{\partial z}[u_r^p(f(z),z)]dz
+\frac{1}{e}h(0)\right |\leq\int _0^{\gamma(p)}|h(z)-h(0)|
\left |\frac{\partial}{\partial z}[u_r^p(f(z),z)]\right |dz$$
$$\mbox{ because }\int _0^{\gamma(p)}\frac{\partial}{\partial z}[u_r^p(f(z),z)]dz=-\frac{1}{e}$$

As $\gamma(p)$ tends to 0 as $p$
tends to infinity, 
$$\forall \varepsilon>0 \exists p_*(\varepsilon)\in\N^*\mbox{ such that }\forall p\geq p_*, \gamma(p)<\delta(\varepsilon)$$ \\
By the same reasoning as the one used to establish
(\ref{divergence}) and as $f'$ is bounded, we are able to show that
$\left |\ds\frac{\partial}{\partial
  z}[u_r^p(f(z),z)]\right |\leq Dp$, with $D$ a real constant.\\
Moreover, as $\ds\gamma(p)\leq\frac{1}{p}$, we deduce that
\begin{multline*}
\forall p\geq p_*(\varepsilon),\\
\left|\int _0^{\gamma(p)}h(z)\ds\frac{\partial}{\partial
  z}[u_r^p(f(z),z)]dz+\frac{1}{e}h(0)\right |\leq\varepsilon\int
_0^{\gamma(p)}\ds\left |\frac{\partial}{\partial
  z}[u_r^p(f(z),z)]\right |dz\\
\leq\varepsilon Dp\ds\frac{1}{p}\leq D\varepsilon\\
\end{multline*}

To summarize 
\begin{multline*}
\forall\varepsilon>0,\exists p_*(\varepsilon)\in\N^*\mbox{ such that }\forall p\geq p_*(\varepsilon),\\
 \left|\int
_0^{\gamma(p)}(h(z)\frac{\partial}{\partial
  z}[u_r^p(f(z),z)]+\frac{1}{e}h(0))dz\right|\leq\varepsilon\\
\end{multline*}
We have also proved that 
\begin{equation}\label{lg_2}
\lim_{p\to +\infty}\int_{\gamma ^*_{LG}}r\textrm{div}
  _{\gamma}(u^p)ds=-\frac{1}{e}h(0)
\end{equation}

\subsubsection{The term $l_{LS}(\omega^*).u^p$}

By (\ref{derivee_bord_ls}), (\ref{div_pol}) and (\ref{scal_pol}), we have 
$l_{LS}(\omega^*).u^p=\mu\ds\int_{\gamma ^*_{LS}}u_r^pds+\mu\ds\int_{\gamma ^*_{LS}}\textrm{div}_{\gamma}(u^p)rdr$

As $\gamma _{LS}\cap\mbox{support}(u^p)=\left [
f(0)-\ds\frac{1}{p},f(0) \right ]$ and $u_r^p$ is bounded by $1$, we have

$$\ds\int_{\gamma ^*_{LS}}u_r^pds\leq(f(0)-(f(0)-\frac{1}{p}))=\frac{1}{p}$$

so
\begin{equation}\label{ls_1}
\lim_{p\rightarrow +\infty}\int_{\gamma ^*_{LS}}u_r^pds=0
\end{equation}

It is the second term containing the surface divergence $\ds\int_{\gamma ^*_{LS}}\textrm{div}_{\gamma}(u^p)rdr$ on $\gamma^*
_{LS}$ which is going to contribute to the calculation of the contact angle, as for the term $l_{LG}$.
To estimate this term, we have to express $^tDu^p$ and $\overrightarrow n$.

We have 
$\ds ^tDu^p=\left ( \begin{array}{cc}
\ds\frac{\partial u^p_r}{\partial r}&0\\
\ds\frac{\partial u^p_r}{\partial z}&0\\
\end{array}\right )$ 
and 
$\ds\overrightarrow n=\left ( \begin{array}{c}
                                         0\\
                                         1\\
                                        \end{array}\right )$\\

so $\ds\int_{\gamma ^*_{LS}}\textrm{div}_{\gamma}(u^p)rdr=\int_{f(0)-\frac{1}{p}}^{f(0)}\frac{\partial u_r^p}{\partial r}(r,0)rdr$\\

Integrating by part, we find
$$\ds\int_{\gamma ^*_{LS}}\textrm{div}_{\gamma}(u^p)rdr=\left
  [ru_r^p(r,0)\right
  ]_{f(0)-\frac{1}{p}}^{f(0)}-
\int_{f(0)-\frac{1}{p}}^{f(0)}u_r^p(r,0)dr
$$

The second term on the right side of this identity has a null limit as $p$ tends to $+\infty$  because $u^p_r$ is bounded by 1.

As $u_r^p(f(0),0)=\frac{1}{e}$ and $u_r^p(f(0)-\frac{1}{p},0)=0$, we have
$$\left [ru_r^p(r,0)\right ]_{f(0)-\frac{1}{p}}^{f(0)}=\frac{f(0)}{e}$$\\

and thus we conclude that 
\begin{equation}\label{ls_2}
\lim_{p\to +\infty}\int_{\gamma
^*_{LS}}\textrm{div}_{\gamma}(u^p)rdr=\frac{f(0)}{e}
\end{equation}

\subsubsection{The term $l_{el}(\omega^*).u^p$}
From (\ref{derivee_el}) and (\ref{div_pol}), we deduce $$l_{el}(\omega^*).u^p=-\frac{\delta}{2}\ds\int_{\omega ^*}\varepsilon u_r^p |\nabla
\varphi^{\omega ^*}|^2drdz-\frac{\delta}{2}\ds\int _{\omega ^*}\varepsilon |\nabla \varphi^{\omega ^*}|^2\textrm{div}(u^p)rdrdz+$$
$$\frac{\delta}{2}\ds\int_{\omega ^*}\varepsilon  <(^tDu^p+Du^p)\nabla \varphi^{\omega ^*};\nabla
\varphi^{\omega ^*}>rdrdz$$

These terms are containing the electrostatic contributions in the derivative of the energy.

For the first term $\ds\int_{\omega ^*}\varepsilon u_r^p |\nabla
\varphi^{\omega ^*}|^2drdz$, we use the dominated convergence's theorem of Lebesgue. Indeed, we have
$$\varepsilon u_r^p\chi_{B_p} |\nabla
              \varphi^{\omega ^*}|^2\rightarrow 0\mbox{ p.p.,}\mbox{
              as }p\rightarrow +\infty $$
Moreover, we can find a positive constant $M$ such that 
$$\varepsilon u_r^p\chi_{B_p} |\nabla
            \varphi^{\omega ^*}|^2\leq Mr|\nabla
              \varphi^{\omega ^*}|^2$$
As we know that $\varphi^{\omega ^*}\in H^1_{\nu}(\omega ^*)$, we deduce that $|\nabla \varphi^{\omega ^*}|^2\in L^1_{\nu}(\omega
              ^*)$.\\

The dominated convergence's theorem of Lebesgue allows us to claim that 
\begin{equation}\label{lcap_1}
\lim_{p\rightarrow +\infty}\int_{\omega ^*}\varepsilon u_r^p |\nabla
\varphi^{\omega ^*}|^2drdz=0
\end{equation}

The second term $\ds\int _{\omega ^*}\varepsilon |\nabla
\varphi^{\omega ^*}|^2\textrm{div}(u^p)rdrdz$ which appears in $l_{cap}$ is more fastidious to treat. It is necessary to compare the behaviour of the two singular terms in the vicinity of $A$ : 
$\textrm{div}(u^p)$ and $|\nabla \varphi^{\omega ^*}|^2$.

We have an upper bound of the divergence of $u^p$ given in
(\ref{divergence}).
For $|\nabla \varphi^{\omega ^*}|^2$, we can make precise the behaviour of its $L^2$ norm in the vicinity of the point $A$ in each domain $\omega ^* _S$ and $\omega ^* _G$.

We choose to work with polar coordinates $(\rho,\theta)$
centered at $A$.
\begin{center}
\scalebox{0.5}{\includegraphics{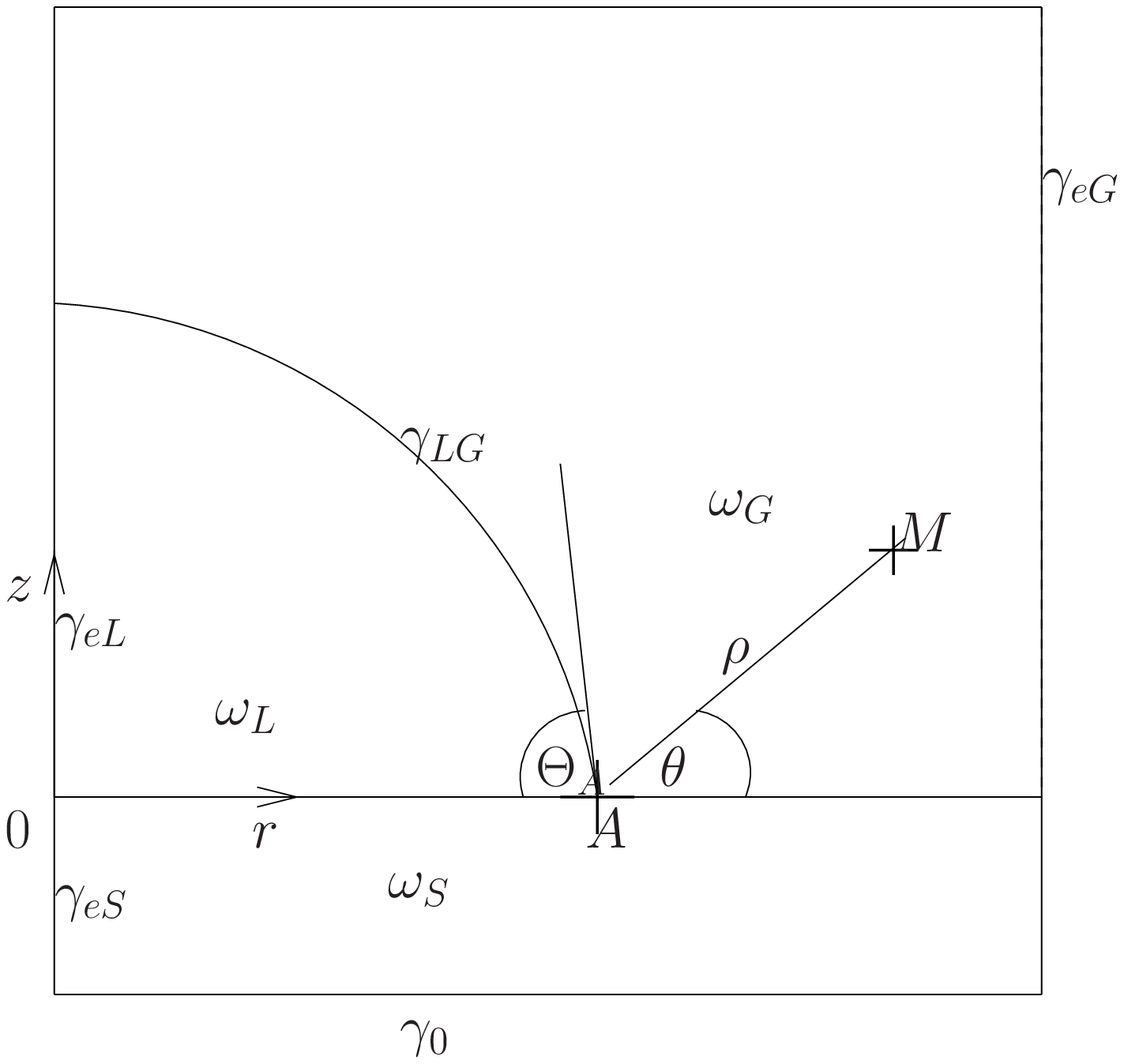}}
\end{center}
Adapting a theorem due to K.Lemrabet (\cite{KL}) which split the
potential into a regular part and a singular part in the vicinity of the edge $A$, we obtain  
\begin{theorem}
It exists a unique $c\in\R$ such that for $i\in\{G,S\}$ : 
$$\varphi^{\omega ^*}_i-c S_i \in H^2_{\nu}(\omega ^* _i)$$
with $$S_G(\rho,\theta)=\frac{\varepsilon _S}{\cos(\lambda(\pi-\Theta
_A))}\rho^{\lambda}\sin(\lambda(\theta -(\pi -\theta
_A)))\eta(\rho,\theta),\quad \theta\in]0,\pi-\Theta _A[$$

$$S_S(\rho,\theta)=\frac{\varepsilon _G}{\cos(\lambda\pi
)}\rho^{\lambda}\sin(\lambda(\theta +\pi
))\eta(\rho,\theta),\quad \theta\in]-\pi,0[$$

where : 

 $\eta$ is an infinitely continuous function with a compact support, such that
$\eta\equiv 1$ in a vicinity of $A$ and becomes zero outside of a ball centered on $A$.

 $\lambda$ is the unique solution in $]0,1[$ of the equation $$\varepsilon
 _S\tan(\lambda(\pi-\Theta _A))=-\varepsilon _G\tan(\lambda\pi)$$.
\end{theorem}

By the equation verified by $\ds\lambda$,
we have $\ds\lambda\in]\frac{1}{2},1[$.

We are now giving an upper bound of the $L^2_{\nu}$ norm of the modulus of the gradient of
$\varphi^{\omega ^*}_i$ on $\omega ^* _i\cap
B_p$ denoted $||\nabla \varphi^{\omega ^*}_i||^2_{L^2_{\nu}(\omega ^*
_i\cap B_p)}$.
By setting $g_i=\varphi^{\omega ^*}_i-c S_i \in H^2_{\nu}(\omega ^* _i)$, we deduce $$||\nabla \varphi^{\omega ^*}_i||^2_{L^2_{\nu}(\omega ^*
_i\cap B_p)}\leq 2\left [||c\nabla S_i||^2_{L^2_{\nu}(\omega ^* _i\cap
B_p)}+||\nabla g_i||^2_{L^2_{\nu}(\omega ^* _i\cap B_p)}\right ]$$
It remains to find an upper bound to the terms which appear in the right hand side of the inequality.

We choose $\tilde p$ big enough such that
$\eta_{/B_{\tilde p}}\equiv
1$ and we consider $p\geq \tilde p$.

Calculating the explicit expression of $\nabla S_i$, we obtain the following upper bound : 
$$\left |\nabla S_i \right |^2\leq D^* \rho ^{2(\lambda -1)}$$
where $D^*$ is a constant.

Then the upper bound of its $L^2_{\nu}$ norm in the vicinity of the point $A$ is
\begin{multline}\label{gradsi}
\left |\left |\nabla S_i \right |\right |^2_{L^2(\omega^*\cap
B_p)}\leq(f(0)+\frac{1}{p})\int_{\omega ^*\cap B_p}\left |\nabla S_i (\rho,\theta) \right
|^2 \rho d\rho d\theta\\
\leq D^*(f(0)+\frac{1}{p})\int_{B_p}\rho ^{2\lambda -1}d\rho d\theta\\
=D^*(f(0)+\frac{1}{p})\frac{\pi}{\lambda}p^{-2\lambda}\\
\end{multline}

Using the Sobolev's injection in three dimension and in axisymmetric case (\cite{SN} p.27)
$$H^1_{\nu}(\omega ^* _i)\subset L^5_{\nu}(\omega ^* _i)$$
we obtain an upper bound of $||\nabla g_i||^2_{L^2_{\nu}(\omega ^* _i\cap
B_p)}$ by Hölder inequality : 
\begin{equation}\label{gradgi}
\int _{\omega ^*_i\cap B_p}|\nabla g_i|^2d\nu=\int_{\omega
^* _i}\underbrace{\chi _{B_p}}_{L^{\frac{5}{3}}_{\nu}(\omega ^* _i)}\underbrace{|\nabla g_i|^2}_{L^{\frac{5}{2}}_{\nu}(\omega ^* _i)}d\nu\leq (f(0)+\frac{1}{p})\underbrace{||\chi
_{B_p}||_{L^{\frac{5}{3}}(\omega ^*
_i)}}_{\mbox{mes}(B_p)^{\frac{3}{5}}}\underbrace{||\nabla
g_i||_{L^5_{\nu}(\omega ^* _i)}^2}_{K_i}
\end{equation}

From the two estimates (\ref{gradsi}) and (\ref{gradgi}), we deduce : 
$$||\nabla \varphi^{\omega ^*}_i||^2_{L^2_{\nu}(\omega ^*
_i\cap B_p)}\leq
2(f(0)+\frac{1}{p})[D^*\frac{\pi}{\lambda}p^{-2\lambda}+K_i
\mbox{mes}(B_p)^{\frac{3}{5}}]$$
for $i\in\{S,G\}$

We can also write an estimation of the gradient norm of the potential on the whole domain $\omega ^*$

\begin{equation}\label{gradphi}
||\nabla \varphi^{\omega ^*}||^2_{L^2_{\nu}(\omega ^*
\cap B_p)}\leq
2(f(0)+\frac{1}{p})\left [D^*\frac{\pi}{\lambda}p^{-2\lambda}+K^*
\mbox{mes}(B_p)^{\frac{3}{5}}\right ]
\end{equation}
where $K^*=K_S+K_G$.

Now we are coming back to the term $\ds\int _{\omega ^*}\varepsilon |\nabla
\varphi^{\omega ^*}|^2\textrm{div}(u^p)rdrdz$.
For every $p\geq \tilde p$, we have
\begin{multline*}
\int _{\omega ^*}\varepsilon |\nabla \varphi^{\omega ^*}|^2\textrm{div}(u^p)rdrdz\leq MDp\underbrace{\int _{\omega ^*\cap B_p}|\nabla
\varphi^{\omega ^*}|^2drdz}_{||\nabla \varphi^{\omega ^*}||^2_{L^2_{\nu}(\omega ^*
\cap B_p)}}\mbox{ par }(\ref{divergence})\\
\leq 4MDp(f(0)+\frac{1}{p})\left [D^*\frac{\pi}{\lambda}p^{-2\lambda}+C^*
\frac{\pi^{\frac{3}{5}}}{p^{\frac{6}{5}}}\right ]\\
\end{multline*}

and as $$-2\lambda+1<0, $$
we deduce that
\begin{equation}\label{lcap_2}
\lim_{p\to +\infty}\int _{\omega ^*}\varepsilon |\nabla
\varphi^{\omega ^*}|^2\textrm{div}(u^p)rdrdz=0
\end{equation}

It remains to examine the term $$\ds\int_{\omega ^*}\varepsilon  <(^tDu^p+Du^p)\nabla \varphi^{\omega ^*};\nabla
\varphi^{\omega ^*}>rdrdz$$ This term is of the same order as the previous term.

We have $\ds^tDu^p+Du^p=\left (\begin{array}{cc}
                         2\frac{\partial u_r^p}{\partial r}&\frac{\partial u_r^p}{\partial z}\\
                         \frac{\partial u_r^p}{\partial z}&0\\
                         \end{array}\right ).$\\

Denoting $\nabla \varphi^{\omega ^*}=\left (\begin{array}{c}
                                   \delta_r\\
                                   \delta_z\\
                                   \end{array}\right )$\\

we obtain

$$<(^tDu^p+Du^p)\nabla \varphi^{\omega ^*};\nabla \varphi^{\omega ^*}>=2\frac{\partial
  u_r^p}{\partial r}(\delta_r)^2+2\frac{\partial u_r^p}{\partial z}\delta_r\delta_z$$

from what we deduce 
\begin{multline*}
\left |\int_{\omega ^*}\varepsilon  <^tDu^p+Du^p)\nabla
\varphi^{\omega ^*};\nabla \varphi^{\omega ^*}>rdrdz\right|\leq 2M\int _{\omega ^*\cap
  B_p}\left ( \left |\frac{\partial u_r^p}{\partial
  r}(\delta_r)^2+\frac{\partial u_r^p}{\partial z}\delta_r\delta_z\right |\right )rdrdz\\
\leq 2M\left [\int _{\omega ^*\cap B_p}|\frac{\partial
  u_r^p}{\partial r}||\nabla \varphi^{\omega ^*}|^2rdrdz+\frac{1}{2}\int
  _{\omega ^*\cap B_p}|\frac{\partial u_r^p}{\partial z}||\nabla
  \varphi^{\omega ^*}|^2rdrdz\right ]\\
\end{multline*}

We have $\ds\left |\frac{\partial u_r^p}{\partial r}\right |\leq 2Cp$, $\ds\left |\frac{\partial u_r^p}{\partial z}\right |\leq Dp$ and $\ds r\leq f(0)+\frac{1}{p}$\\

By an analogous reasoning, we can show that 
\begin{equation}\label{lcap_3}
\lim_{p\to +\infty}\int_{\omega ^*}\varepsilon  <^tDu^p+Du^p)\nabla
\varphi^{\omega ^*};\nabla \varphi^{\omega ^*}>rdrdz=0
\end{equation}

\subsubsection{Conclusion}

In the same manner, we can prove that 
\begin{equation}\label{lcont}
\lim _{p\to +\infty}l_{cont}(\omega ^*).u^p=0
\end{equation}
From (\ref{grav_3})-(\ref{lcap_1}), (\ref{lcap_2})-(\ref{lcont}), we deduce 
\begin{equation}\label{limite}
\lim_{p \to +\infty}l(\omega ^*).u^p=\frac{1}{e}[-h(0)+\mu f(0)]
\end{equation}

\subsection{Value of the wetting angle}

At the minimum of energy $\omega ^*$, we verify : 
$$\exists \lambda ^*\in\R \mbox{ such that for all }p\mbox{ big enough},\quad l(\omega ^*).u^p=\lambda ^*l_{cont}(\omega^*).u^p.$$
So 
$$\lim_{p\to +\infty}l(\omega^*).u^p=\lambda ^*\lim _{p\to +\infty}l_{cont}(\omega^*).u^p$$
From (\ref{limite}) we deduce that 
$$h(0)=\mu f(0)$$
And so if $f(0)\neq 0$, we conclude : 
$$\mu=\frac{f'(0)}{\sqrt{1+f'^2(0)}}$$

But we know that the right hand side is in fact the cosine of the contact angle

$$\cos(\Theta _A)=-\mu$$

By the definition of $\mu$ and Young's angle, we deduce that the value of the contact angle for an optimal drop with an axisymmetric geometry is Young's angle and this is valid for all values of potential applied.

\section{Numerical work in progress}\label{num_work}
 We recall that a simple model consists in considering that the system is a plane capacitor. In this case, the contact angle is given by the relation : 
$$\cos(\theta)=\cos(\theta_Y)+\frac{\varepsilon_S}{2e\sigma_{LG}}V^2$$
where $e$ is the thickness of the insulator. It predicts a total spreading of the drop on the polymer. In fact we can admit that it gives a value of a macroscopic contact angle.

We use two numerical approaches. The first one is a macroscopic one and the second one a local one. 

From the code "Electrocap" developed by J.Monnier and P.Chow Wing Bom (\cite{MW}, \cite{MW_bis}), we introduce a treatment of the singularity of the potential.
In order to compute the singularity of the potential we use the Singular Complement Method as presented in \cite{Ciarlet}.

This numerical approach gives same qualitative results as in \cite{MW_bis} (which is without the treatment of the singularity). It shows a deviation from the shape predicted by the plane capacitor approximation : the contact angle is higher than predicted by the plane capacitor approximation and the curvature increases sharply near the triple line. But the effect of the treatment of the singularity seems to be deleted. Furthermore it doesn't show that the contact angle is constant.

We present numerical results for given volume and physical parameters.

\begin{figure}[htbp]
\begin{tabular}{cc}
\includegraphics[width=0.5\linewidth]{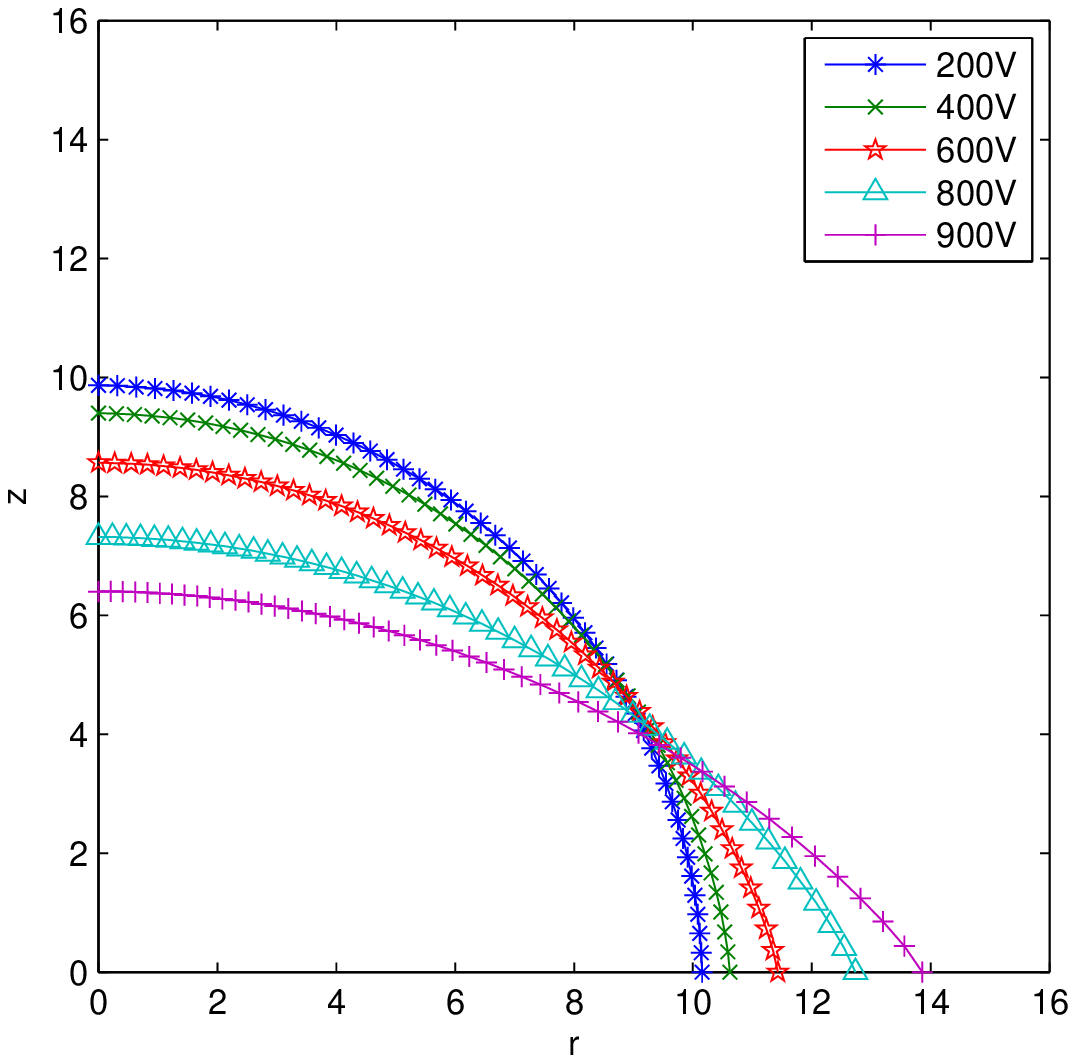} & 
\includegraphics[width=0.5\linewidth]{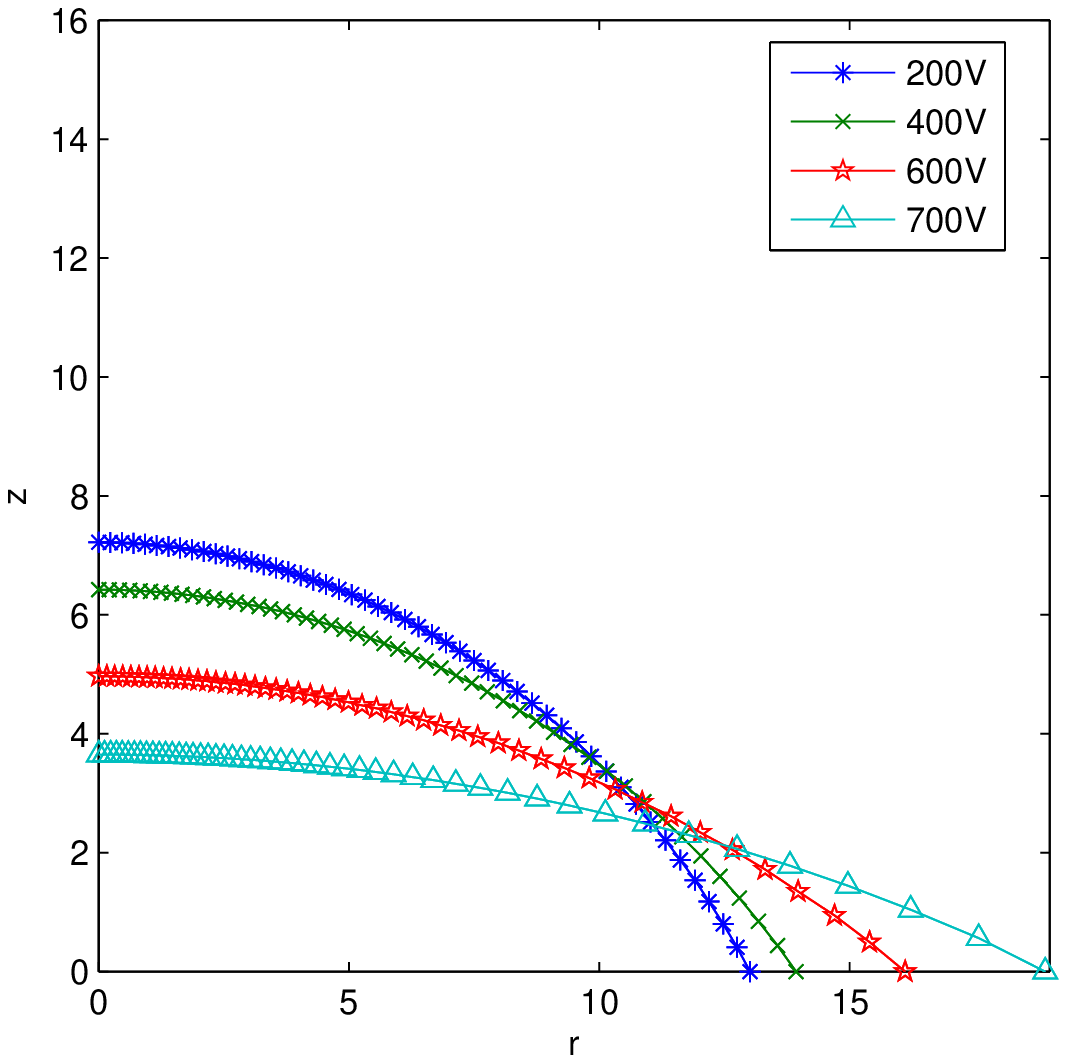}
\end{tabular}
\caption{{\it Left, }Macroscopic shape of the drop for a Young's angle of $90°$;{\it Right, }Macroscopic shape of the drop for a Young's angle of $60°$}
\label{shapes}
\end{figure}

In Fig.\ref{shapes} we present the macroscopic shape of the drop obtained for different voltages.

A second numerical work is also in progress. The global approximation (using the treatment of the singularity) is used together with a local model (given by an ODE in the vicinity of the contact point).

\begin{figure}[htbp]
\begin{tabular}{cc}
\includegraphics[width=0.5\linewidth]{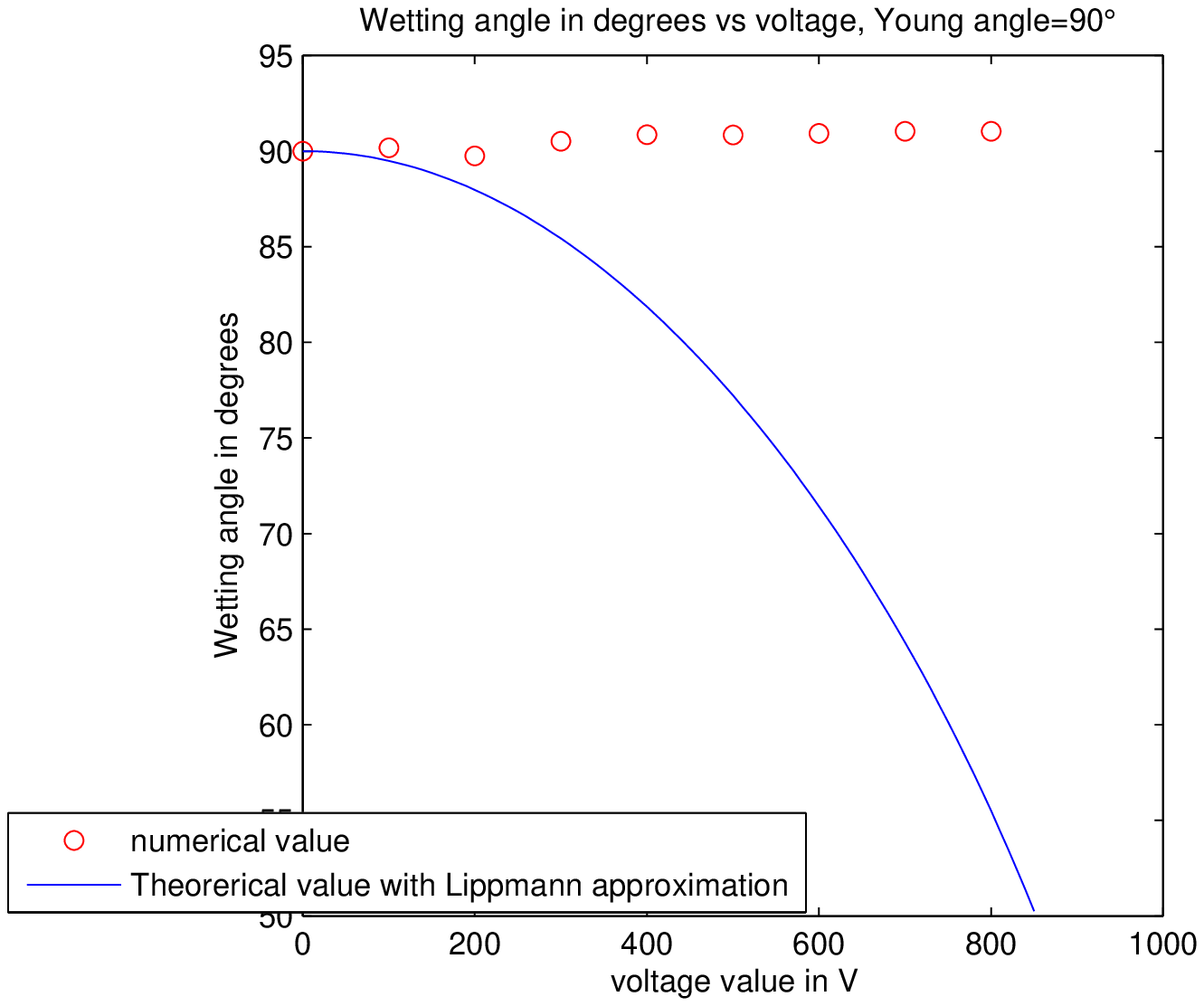} & 
\includegraphics[width=0.5\linewidth]{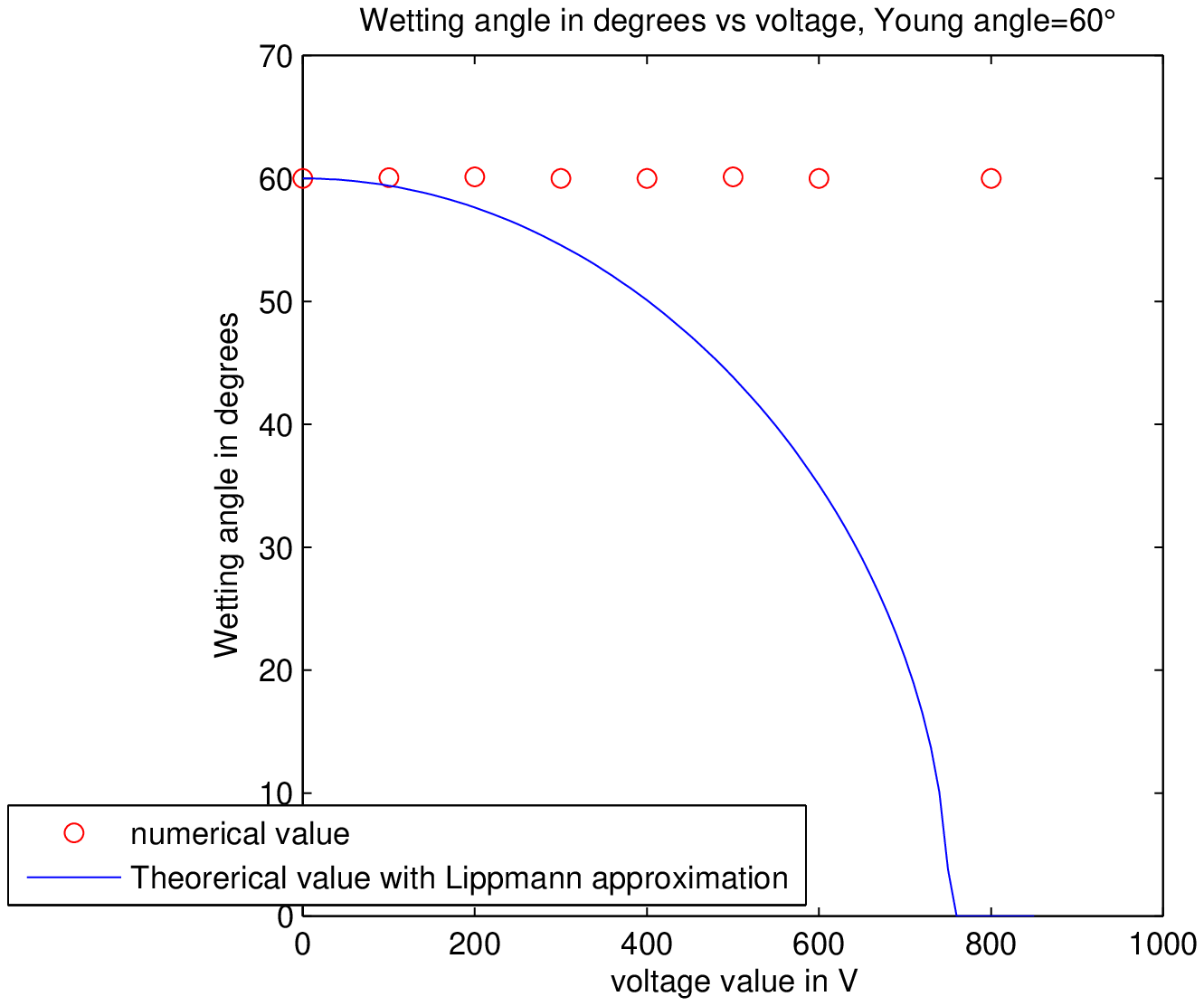}
\end{tabular}
\caption{{\it Left, }Value of the contact angle for a Young angle of $90°$; {\it Right, }Value of the contact angle for a young angle of $60°$}
\label{contact_angles}
\end{figure}

In Fig.\ref{contact_angles} are presented values of the contact angle obtained by this numerical approach. Values are compared to the plane capacitor approximation.\\

\begin{figure}[htbp]
\begin{tabular}{cc}
\includegraphics[width=0.5\linewidth]{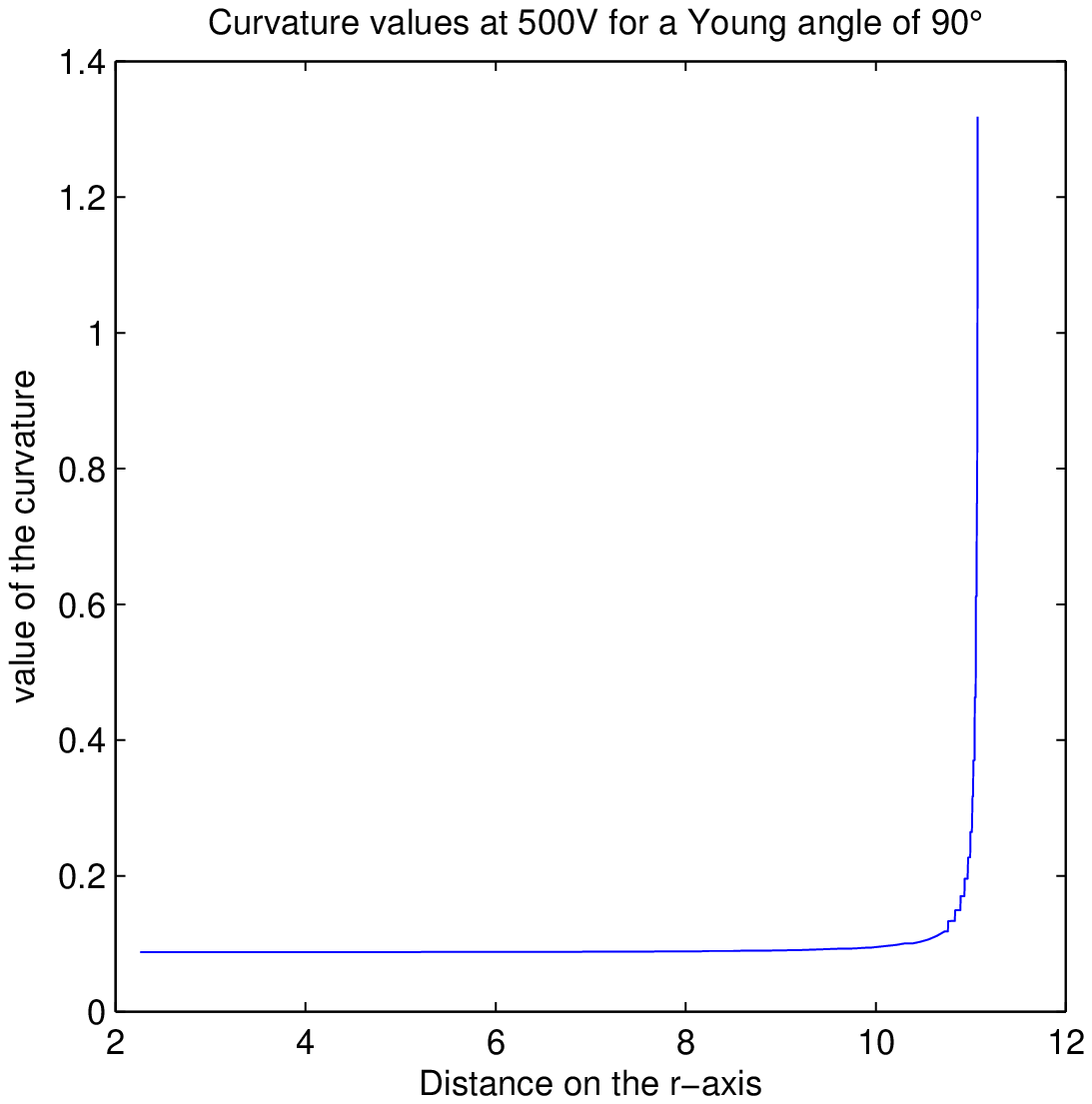} & 
\includegraphics[width=0.5\linewidth]{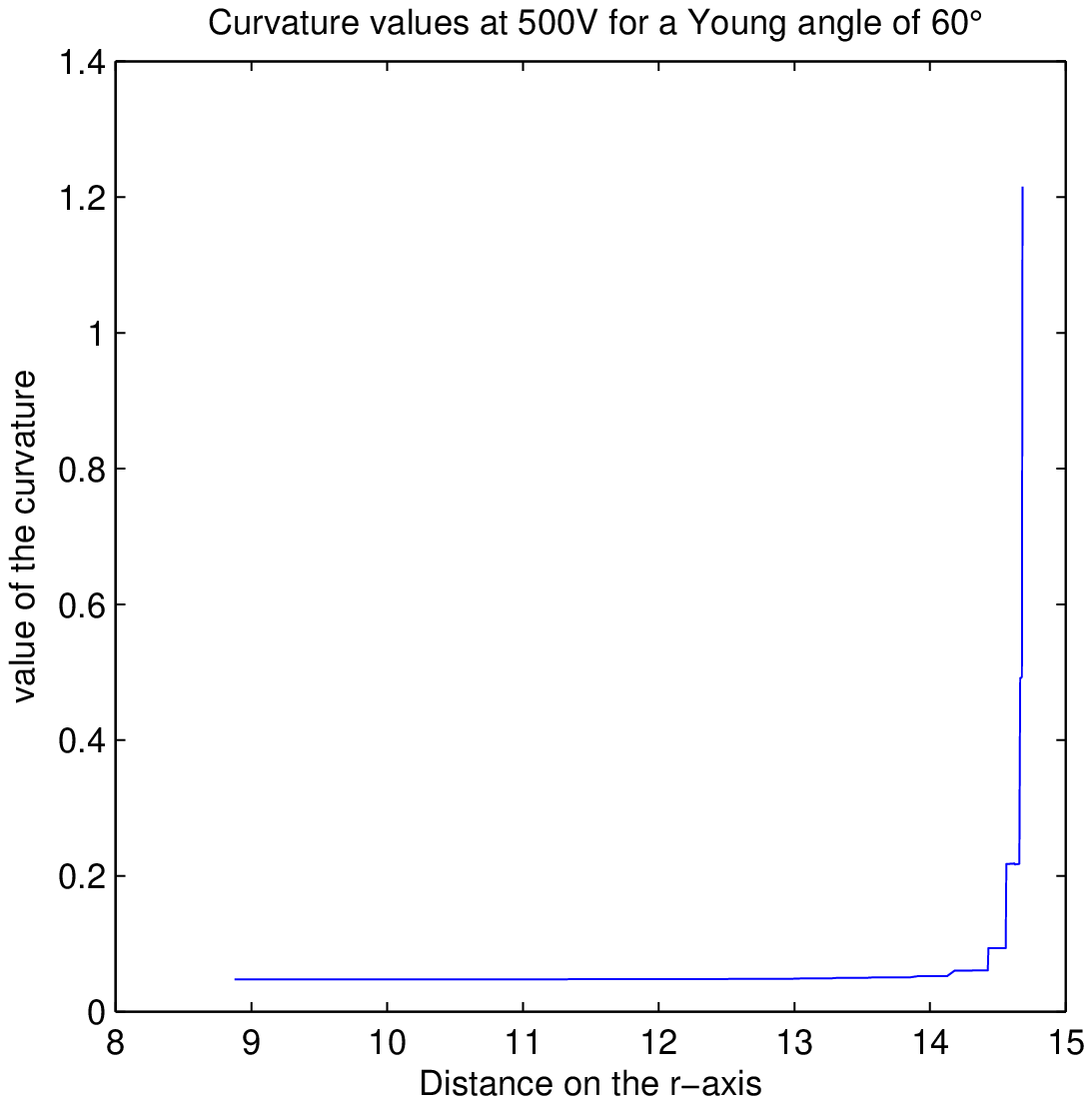}
\end{tabular}
\caption{{\it Left, }Value of the curvature for a Young angle of $90°$ at 500V; {\it Right, }Value of the curvature for a Young angle of $60°$ at 500V}
\label{curvature}
\end{figure}

The numerical contact angle appears to be constant as theoretically predicted in the last section and proposed in \cite{allemands}.
The numerical curvature appears to explode near the contact point. In Fig.\ref{curvature} we present the curvature value for points of the drop at $500V$ for a Young's angle of $90°$. For other voltages and values of Young's angle, the qualitative behaviour remains the same.\\

All this results are in accordance with the theory and experimental works (see e.g.\cite{MB}).

 \section{Conclusion}\label{conclusion}
 We have proved that the contact angle in electrowetting remains constant for all potential applied and it equals Young's angle. It is the result that has been predicted in \cite{allemands}. A numerical work is also in progress. Because of the singularity of the potential in the vicinity of the contact point, it uses a local model near the triple line. The contact angle computed appears to be constant.

 \end{document}